\documentclass[12pt,reqno]{amsart}
\usepackage{a4wide}
\usepackage{amsfonts,amsmath,amssymb}

\title{Human proofs of identities by Osburn and Schneider}

\author[H.~Prodinger]{Helmut Prodinger}
\address{Helmut Prodinger\\
Department of Mathematics\\
University of Stellenbosch\\
7602 Stellenbosch\\
South Africa}
\email{hproding@sun.ac.za}

\date{\today}
\begin{document}
\begin{abstract}
Osburn and Schneider~\cite{OsSc06} derived several combinatorial identities involving
harmonic numbers using the computer programm \textsf{Sigma}. Here, they are derived
by partial fraction decomposition and creative telescoping.
\end{abstract}

\maketitle

\section{Introduction}

In \cite{OsSc06}, sums of the type
\begin{equation*}
\sum_{k=0}^n\binom nk\binom{n+k}k(-1)^{n-k}\cdot X
\end{equation*}
where computed, where $X$ was one of the list
\begin{equation*}
1,\ H_k,\ kH_k,\ kH_{n+k},\ kH_{n-k},\ H_k^{(2)},\ldots,
\end{equation*}
with harmonic numbers of first and second order. The powerful computer software
\textsf{Sigma} was used to evaluate these sums.

In this note, I want to show a method that humans can use; in this way we obtain
``computer-free'' proofs (although it is always helpful to use a computer, to ease computations). The method
consists of the following steps:

\begin{itemize}
\item A rational function, of the form
\begin{equation*}
\frac{(z+1)\dots(z+n)}{z(z-1)\dots(z-n)}\cdot Y
\end{equation*}
will undergo \emph{partial fraction decomposition}.

\item The resulting equation will be multiplied by $z$, and the limit $z\to\infty$
will be performed. In the simpler cases, this is already the desired identity.

\item If, say, $H_k$ is involved, we write it as
\begin{equation*}
H_k=\sum_{j\ge1}\bigg[\frac1j-\frac1{j+k}\bigg]=\sum_{j\ge1}\frac{k}{j(j+k)}.
\end{equation*}
We use an individual term $\frac{z}{j(j+z)}$ as part of the rational function. The resulting
identities will be summed over all $j\ge1$.

\item To obtain the final result, these infinite series have to be evaluated. One possible
route would be to use hypergeometric series and the huge collection of identities and transformations that is known~\cite{Krattenthaler93}; the desired formula usually comes out as a limiting case.

However, here, we decided to use \emph{creative telescoping}, as popularized in the book
\cite{PeWiZe96}.
\end{itemize}

This method was partially inspired by \cite{Chu04}.

I don't want to go through the whole list of identities by Osburn and Schneider, but a fair amount
of them of increasing difficulty will be treated, to convince the reader, that there are no principal
limitations to the approach. 

\section{Identities}

\subsection{} 

Consider
\begin{equation*}
\frac{(z+1)\dots(z+n)}{z(z-1)\dots(z-n)}=\sum_{k=0}^n\binom nk\binom{n+k}k(-1)^{n-k}
\frac1{z-k}.
\end{equation*}
Multiplying this by $z$, and letting $z\to\infty$, we obtain
\begin{equation*}
\sum_{k=0}^n\binom nk\binom{n+k}k(-1)^{n-k}=1.
\end{equation*}

\subsection{} 

Consider
\begin{multline*}
\frac{(z+1)\dots(z+n-1)}{z(z-1)\dots(z-n)}\frac1{z+n}
\\=
\sum_{k=0}^n\binom nk\binom{n+k}k(-1)^{n-k}\frac1{(n+k)^2}\frac1{z-k}
+\frac{(n-1)!^2}{(2n)!}\frac{1}{z+n}.
\end{multline*}
The limit form is
\begin{equation*}
\sum_{k=0}^n\binom nk\binom{n+k}k(-1)^{n-k}\frac1{(n+k)^2}
=-\frac{(n-1)!^2}{(2n)!}.
\end{equation*}

\subsection{} 

\begin{multline*}
\frac{(z+1)\dots(z+n)}{(z-1)\dots(z-n)}\frac{1}{j(j+z)}\\
=\sum_{k=1}^n\binom nk\binom{n+k}k(-1)^{n-k}\frac{k}{j(j+k)}\frac1{z-k}+
\frac{(j-1)!^2}{(j-n-1)!(n+j)!}\frac{1}{j+z}.
\end{multline*}
The limit form is
\begin{align*}
\sum_{k=1}^n\binom nk\binom{n+k}k(-1)^{n-k}\frac{k}{j+k}=
-\frac{(j-1)!^2}{(j-n-1)!(n+j)!}+\frac{1}{j}.
\end{align*}
Summing on $j\ge1$ (and shifting the index), we get
\begin{align*}
\sum_{k=1}^n\binom nk\binom{n+k}k(-1)^{n-k}H_k&=
\sum_{j\ge0}\bigg[\frac{1}{j+1}-\frac{j!^2}{(j-n)!(n+1+j)!}\bigg].
\end{align*}
Denote the summand of the last expression by $F(n,j)$. Then we have
\begin{equation*}
F(n+1,j)-F(n,j)=G(n,j+1)-G(n,j),
\end{equation*}
with
\begin{equation*}
G(n,j)=\frac{2j!^2}{(j-n-1)!(n+1+j)!(n+1)}.
\end{equation*}
Now set $S(n):=\sum_{j\ge0}F(n,j)$, then
\begin{equation*}
S(n+1)-S(n)=\lim_{J\to\infty}G(n,J)-G(n,0)=\frac2{n+1}.
\end{equation*}
Together with $S(0)=0$, this leads to $S(n)=2H_n$, and the identity
\begin{equation*}
\sum_{k=1}^n\binom nk\binom{n+k}k(-1)^{n-k}H_k=2H_n.
\end{equation*}

\subsection{} 

Consider
\begin{equation*}
A:=\frac{(z+1)\dots(z+n)}{z\dots(z-n)}\frac{n+z}{j(j+n+z)}.
\end{equation*}
Then
\begin{align*}
\sum_{k=0}^n\binom nk\binom{n+k}k(-1)^{n-k}\frac{n+k}{j(j+n+k)}\frac1{z-k}=A-
\frac{(n+j-1)!^2}{(j-1)!(2n+j)!}\frac{1}{j+n+z}
\end{align*}
and
\begin{align*}
\sum_{k=0}^n\binom nk\binom{n+k}k(-1)^{n-k}\frac{n+k}{j(j+n+k)}=\frac1j-
\frac{(n+j-1)!^2}{(j-1)!(2n+j)!}.
\end{align*}
Summing,
\begin{align*}
\sum_{k=0}^n\binom nk\binom{n+k}k(-1)^{n-k}H_{n+k}
&=\sum_{j\ge0}\bigg[\frac1{j+1}-\frac{(n+j)!^2}{j!(2n+1+j)!}\bigg].
\end{align*}
Denoting the term in the sum again by $F(n,j)$, we find
\begin{equation*}
F(n+1,j)-F(n,j)=G(n,j+1)-G(n,j),
\end{equation*}
with
\begin{equation*}
G(n,j)=\frac{(n+j)!^2(3n+3+2j)}{(j-1)!(2n+2+j)!(n+1)}.
\end{equation*}
For the summatory functions, we thus get
\begin{equation*}
S(n+1)-S(n)=\frac{2}{n+1},
\end{equation*}
with the same result as before, and henceforth the identity
\begin{equation*}
\sum_{k=1}^n\binom nk\binom{n+k}k(-1)^{n-k}H_{n+k}=2H_n.
\end{equation*}

\subsection{} 
Consider
\begin{multline*}
\frac{(z+1)\dots(z+n)}{(z-1)\dots(z-n)}\frac{z}{j(j+z)}\\
=\sum_{k=1}^n\binom nk\binom{n+k}k(-1)^{n-k}\frac{k^2}{j(j+k)}\frac1{z-k}-\frac{(j-1)!j!}{(j-n-1)!(n+j)!}\frac{1}{j+z}+\frac1j
\end{multline*}
and the limiting form
\begin{align*}
\sum_{k=1}^n\binom nk\binom{n+k}k(-1)^{n-k}\frac{k^2}{j(j+k)}=\frac{n(n+1)}{j}-1+\frac{(j-1)!j!}{(j-n-1)!(n+j)!}.
\end{align*}
Summing on $j$:
\begin{align*}
\sum_{k=1}^n\binom nk\binom{n+k}k(-1)^{n-k}kH_k
=\sum_{j\ge0}\bigg[\frac{n(n+1)}{j+1}-1+\frac{j!(j+1)!}{(j-n)!(n+1+j)!}\bigg].
\end{align*}
Hence
\begin{equation*}
(n+2)F(n,j)-nF(n+1,j)=G(n,j+1)-G(n,j)-2,
\end{equation*}
with
\begin{equation*}
G(n,j)=\frac{2j!(j+1)!}{(j-n-1)!(n+1+j)!}.
\end{equation*}
Thus
\begin{equation*}
(n+2)S(n)-nS(n+1)=\lim_{J\to\infty}G(n,J)-2J-G(n,0)=-2n(n+2),
\end{equation*}
and furthermore
\begin{equation*}
\frac{S(n+1)}{(n+1)(n+2)}-\frac{S(n)}{n(n+1)}=\frac{2}{n+1}.
\end{equation*}
Unwinding this, we get
\begin{equation*}
\frac{S(n)}{n(n+1)}=\frac{S(1)}{2}+\sum_{k=1}^{n-1}\frac{2}{k+1}=1+2H_n-2.
\end{equation*}
So we found $S(n)=n(n+1)(2H_n-1)$, and the identity
\begin{equation*}
\sum_{k=1}^n\binom nk\binom{n+k}k(-1)^{n-k}kH_k
=n(n+1)(2H_n-1).
\end{equation*}

\subsection{}

Consider
\begin{multline*}
\frac{(z+1)\dots(z+n)}{(z-1)\dots(z-n)}\frac{n+z}{j(j+n+z)}\\
=\sum_{k=1}^n\binom nk\binom{n+k}k(-1)^{n-k}\frac{k(n+k)}{j(j+n+k)}\frac1{z-k}-\frac{(n+j-1)!(n+j)!}{(j-1)!(2n+j)!}\frac{1}{j+n+z}+\frac1j
\end{multline*}
and the limiting form
\begin{align*}
\sum_{k=1}^n\binom nk\binom{n+k}k(-1)^{n-k}\frac{k(n+k)}{j(j+n+k)}=\frac{n(n+1)}{j}-1+
\frac{(n+j-1)!(n+j)!}{(j-1)!(2n+j)!}
\end{align*}
which we sum:
\begin{align*}
\sum_{k=1}^n\binom nk\binom{n+k}k(-1)^{n-k}kH_{n+k}
=\sum_{j\ge0}\bigg[\frac{n(n+1)}{j+1}-1+\frac{(n+j)!(n+1+j)!}{j!(2n+1+j)!}\bigg].
\end{align*}
Hence
\begin{equation*}
(n+2)F(n,j)-nF(n+1,j)=G(n,j+1)-G(n,j)-2,
\end{equation*}
with
\begin{equation*}
G(n,j)=\frac{(n+j)!(n+1+j)!(3n+4+2j)}{(j-1)!(2n+2+j)!}.
\end{equation*}
Consequently
\begin{equation*}
(n+2)S(n)-nS(n+1)=\lim_{J\to\infty}G(n,J)-2J-G(n,0)=-n(2n+3).
\end{equation*}
The solution of this is $S(n)=2n(n+1)H_n-n^2$, and thus we have the identity
\begin{equation*}
\sum_{k=1}^n\binom nk\binom{n+k}k(-1)^{n-k}kH_{n+k}=2n(n+1)H_n-n^2.
\end{equation*}

\subsection{}
Consider
\begin{multline*}
\frac{(z+1)\dots(z+n)}{(z-1)\dots(z-n)}\frac{n-z}{j(j+n-z)}\\
=\sum_{k=1}^n\binom nk\binom{n+k}k(-1)^{n-k}\frac{k(n-k)}{j(j+n-k)}\frac1{z-k}+\frac{(2n+j)!(j-1)!}{(n+j-1)!(n+j)!}\frac{1}{j+n-z}+\frac1j
\end{multline*}
and the limiting form
\begin{align*}
\sum_{k=1}^n\binom nk\binom{n+k}k(-1)^{n-k}\frac{k(n-k)}{j(j+n-k)}=\frac{n(n+1)}{j}+1-
\frac{(2n+j)!(j-1)!}{(n+j-1)!(n+j)!},
\end{align*}
which we sum:	
\begin{align*}
\sum_{k=1}^n\binom nk\binom{n+k}k(-1)^{n-k}kH_{n-k}=\sum_{j\ge0}\bigg[\frac{n(n+1)}{j+1}+1-\frac{(2n+j+1)!j!}{(n+j)!(n+j+1)!}\bigg].
\end{align*}
We have
\begin{equation*}
-nF(n+1,j)+(n+2)F(n,j)=G(n,j+1)-G(n,j)+2,
\end{equation*}
with
\begin{equation*}
G(n,j)=-\frac{(2n+j+1)!j!(3n+2+2j)}{(n+j)!(n+j+1)!}.
\end{equation*}
Hence
\begin{align*}
-nS(n+1)+(n+2)S(n)&=\lim_{J\to\infty}G(n,J)+2J-G(n,0)\\
&=-(n+2)(2n+1)+\frac{(2n+1)!}{n!(n+1)!}(3n+2),
\end{align*}
or
\begin{align*}
-\frac{S(n+1)}{(n+1)(n+2)}+\frac{S(n)}{n(n+1)}
&=-\frac{2n+1}{n(n+1)}+\frac{(2n+1)!}{(n+1)!(n+2)!n}(3n+2),
\end{align*}
and
\begin{align*}
-\frac{S(n)}{n(n+1)}+\frac{S(1)}{2}
&=-\sum_{k=1}^{n-1}\frac{2k+1}{k(k+1)}+\sum_{k=1}^{n-1}\frac{(2k+1)!}{(k+1)!(k+2)!k}(3k+2).
\end{align*}
Thus
\begin{align*}
S(n)&=2n(n+1)H_n-(n+1)^2-n(n+1)\sum_{k=1}^{n-1}\frac{(2k+1)!}{(k+1)!(k+2)!k}(3k+2).
\end{align*}
So we derived the identity
\begin{multline*}
\sum_{k=1}^n\binom nk\binom{n+k}k(-1)^{n-k}kH_{n-k}\\=
2n(n+1)H_n-(n+1)^2-n(n+1)\sum_{k=1}^{n-1}\frac{(2k+1)!}{(k+1)!(k+2)!k}(3k+2).
\end{multline*}
The paper \cite{OsSc06} gives the different evaluation
\begin{align*}
2n(n+1)H_n-(n+1)^2+(2n+1)\binom{2n}n-\frac32n(n+1)\sum_{k=1}^{n}\frac{(2k)!}{k!(k+1)!},
\end{align*}
but it is a routine check to prove that they are the same.

\subsection{}

Consider
\begin{align*}
\frac{(z+1)\dots(z+n)}{(z-1)\dots(z-n)}&\frac{2j+z}{j^2(j+z)^2}=
\sum_{k=1}^n\binom nk\binom{n+k}k(-1)^{n-k}\frac{k(k+2j)}{j^2(j+k)^2}\frac1{z-k}
\\&+\frac{(j-1)!^2}{(n+j)!(j-n-1)!}\frac1{(j+z)^2}\\&+
\frac{(j-1)!^2}{(n+j)!(j-n-1)!}(H_{j+n}+H_{j-n-1}-2H_{j-1})\frac1{j+z}
\end{align*}
and the limiting form
\begin{align*}
\frac1{j^2}&=
\sum_{k=1}^n\binom nk\binom{n+k}k(-1)^{n-k}\frac{k(k+2j)}{j^2(j+k)^2}
+\frac{(j-1)!^2}{(n+j)!(j-n-1)!}(H_{j+n}+H_{j-n-1}-2H_{j-1}).
\end{align*}
Summing,
\begin{align*}
\sum_{k=1}^n\binom nk\binom{n+k}k(-1)^{n-k}H_k^{(2)}&
=\sum_{j\ge0}\bigg[\frac1{(j+1)^2}-\frac{j!^2}{(n+1+j)!(j-n)!}(H_{j+1+n}+H_{j-n}-2H_{j})\bigg]\\
&=\sum_{j\ge0}\bigg[\frac1{(j+1)^2}+\frac d{dx}\frac{(j+x)!^2}{(n+1+j+x)!(j+x-n)!}\bigg|_{x=0}\bigg].
\end{align*}
Now set
\begin{equation*}
F(n,j):=\frac x{(j+1)^2}+\frac{(j+x)!^2}{(n+1+j+x)!(j+x-n)!},
\end{equation*}
and
\begin{equation*}
G(n,j)=\frac{2(j+x)!^2}{(n+1+j+x)!(j+x-n-1)!},
\end{equation*}
then
\begin{equation*}
-(n+1)F(n,j)+(n+1)F(n+1,j)=G(n,j+1)-G(n,j).
\end{equation*}
Summing,
\begin{equation*}
-(n+1)S(n)+(n+1)S(n+1)=\lim_{J\to\infty}G(n,J)-G(n,0)=2-\frac{2x!^2}{(n+1+x)!(x-n-1)!}.
\end{equation*}
Further,
\begin{equation*}
S(n)-S(n-1)=\frac2n-\frac{2x!^2}{(n+x)!(x-n)!n}.
\end{equation*}
Now we differentiate this relation w.r.t. $x$ and set $x=0$. With an obvious notation, this yields
\begin{equation*}
T(n)-T(n-1)=\frac{2(-1)^{n-1}}{n^2},
\end{equation*}
and
\begin{equation*}
T(n)=2\sum_{k=1}^n\frac{(-1)^{k-1}}{k^2}.
\end{equation*}
Hence we have proved the identity
\begin{equation*}
\sum_{k=1}^n\binom nk\binom{n+k}k(-1)^{n-k}H_k^{(2)}
=2\sum_{k=1}^n\frac{(-1)^{k-1}}{k^2}.
\end{equation*}

	\bibliographystyle{plain} 

%\bibliography{carsten}

\end{document}